\documentclass[11pt]{article} \usepackage{amssymb}
\usepackage{amsfonts} \usepackage{amsmath} \usepackage{amsthm} \usepackage{bm}
\usepackage{latexsym} \usepackage{epsfig}
\usepackage{algorithm}
\usepackage{algorithmic}
\usepackage{graphicx}
\usepackage{url}

\newtheorem*{theorem*}{Theorem}
\newtheorem{theorem}{Theorem}[section]
\newtheorem{lemma}[theorem]{Lemma}
\newtheorem*{proposition*}{Proposition}

\newtheorem{definition}[theorem]{Definition}

\newtheorem{example}[theorem]{Example}
\newcommand{\ignore}[1]{}

\newcommand{\enote}[1]{} \newcommand{\knote}[1]{}
\newcommand{\rnote}[1]{}

\renewcommand{\phi}{\varphi}

\begin{document}
\title{Complete Characterization of Functions Satisfying the Conditions of
Arrow's Theorem}
\author{
  Elchanan Mossel\footnote{Weizmann Institute and U.C. Berkeley. E-mail:
    mossel@stat.berkeley.edu. Supported by a Sloan fellowship in
    Mathematics, by BSF grant 2004105, by NSF Career Award (DMS 054829)
    by ONR award N00014-07-1-0506 and by ISF grant 1300/08}~ and
  Omer Tamuz\footnote{Weizmann Institute. Supported by ISF grant 1300/08}
 }

\date{\today}
\maketitle
\begin{abstract}
  Arrow's theorem implies that a social welfare function satisfying
  Transitivity, the Weak Pareto Principle (Unanimity) and Independence
  of Irrelevant Alternatives (IIA) must be dictatorial.  When
  non-strict preferences are also allowed, a dictatorial social
  welfare function is defined as a function for which there exists a
  single voter whose strict preferences are followed.  This definition
  allows for many different dictatorial functions, since non-strict
  preferences of the dictator are not necessarily followed.  In
  particular, we construct examples of dictatorial functions which do
  not satisfy Transitivity and IIA. Thus Arrow's theorem, in the case
  of non-strict preferences, does not provide a complete
  characterization of all social welfare functions satisfying
  Transitivity, the Weak Pareto Principle, and IIA.

  The main results of this article provide such a characterization for
  Arrow's theorem, as well as for follow up results by Wilson.  In
  particular, we strengthen Arrow's and Wilson's result by giving an
  exact if and only if condition for a function to satisfy
  Transitivity and IIA (and the Weak Pareto Principle). Additionally,
  we derive formulae for the number of functions satisfying these
  conditions.

\end{abstract}

\section{Introduction}

\subsection{Arrow's Model and Non-Strict Preferences}
Arrow's famous impossibility theorem~\cite{Arrow50,Arrow63} shows that
a group of voters, each with a preference regarding a set of
alternatives, has a limited number of ways to come to an agreement. In
fact, under some restrictions, the only possible system is {\em
  dictatorship}, that is, the unqualified adoption of one particular
voter's preference.

In Arrow's model, a set of voters each have a preference, or ranking,
of a set of social alternatives. A {\em Social Welfare Function} (SWF)
is a function from a set of voters' preferences to a single,
``agreed'' preference. The preferences (both of the voters and of the
outcome) may not be strict, so two alternatives could be equally
ranked. However, as may be expected, each preference relation is {\em
  transitive} in the usual sense that if alternative $x$ is preferable
to $y$, and $y$ is preferable to $z$, then $x$ is preferable to $z$.

By Arrow's definition, there exists in a {\em dictatorship} a single
voter whose preferences are always adopted by the SWF. However, the
SWF is bound only to the dictator's {\bf strict} preferences, so that
if the dictator decrees that $x$ is equivalent to $y$, then the SWF
may choose any of the three possible relations between $x$ and $y$.

Arrow's theorem states that an SWF whose outcome is transitive, and
which also satisfies IIA and the Weak Pareto Principle, must be
dictatorial.  Recall that {\em Independence of Irrelevant
  Alternatives} (IIA) requires that the relative ranking of two
alternatives, as decided by the SWF, depend only on the relative
rankings assigned to these alternatives by the voters. Hence, if the
SWF, given a set of voters' preferences, prefers alternative $x$ to
$y$, then it must also prefer $x$ to $y$ for an altered set of
preferences, given that the relative rankings of $x$ and $y$ remains
unchanged.

The second condition, the {\em Weak Pareto Principle} (WPP), asserts
that when all voters have the same strict preference, then the SWF
must concur and the agreed ranking equals this common ranking.

\subsection{What Happens when the Dictator Can't Decide?}
If one allows strict rankings only, then all dictatorial functions
satisfy IIA and WPP, and so Arrow's theorem provides a complete ``if
and only if'' characterization: an SWF satisfies IIA and WPP iff it is
dictatorial. However, in the general case of non-strict preferences
``{\em If the dictator is indifferent between {\em x} and {\em y},
  presumably he will then leave the choice up to some or all of the
  other members of society,}'' in the words of
Arrow~\cite{Arrow50}. Hence, many different dictatorial functions are
possible, and in particular ones that don't satisfy IIA (we provide
examples in Sect.~\ref{section:theorem} below). We are interested in
finding a full characterization of the SWFs satisfying Arrow's
theorem.

\subsection{Wilson's Framework}
We approach this problem under the slightly more general framework of
Wilson~\cite{Wilson72}, who proves an extension of Arrow's
theorem. Further weakening the Weak Pareto Principle, Wilson assumes
``Citizens' Sovereignty'': that there exists {\bf no} pair of
alternatives $x$ and $y$ for which the SWF always deems $x$ strictly
preferable to $y$, regardless of the voters' preferences. Under this
weaker restriction, Wilson shows, in his first theorem, that a
dictatorship still arises, albeit with a slightly modified definition
of a dictator. Under Wilson's conditions, there exists one voter whose
preferences are either always followed verbatim or always followed
{\em in reverse}.

In his second theorem, Wilson disposes with the second condition
entirely, leaving IIA only. He characterizes SWFs which satisfy IIA,
but, like Arrow, does not give a full ``if and only if'' condition. In
particular, and again like Arrow, he does not elaborate on how an SWF
may decide on pairs which a dictator deems equal.

\subsection{Full Characterization}
We complete these gaps and give a full characterization of SWFs that
satisfy the IIA condition. We show that functions satisfying Arrow's
conditions defer decision, on pairs which the dictator deems
equivalent, to a particular class of functions which we dub
``clerical dictatorial''. These functions divide each set of
equivalent alternatives to an ordered partition, let a dictator decide
on each partition and recursively defer the dictators' undecided
pairs to other ``clerical dictatorial'' functions. Eventually,
undecided pairs which are otherwise ordered with respect to the other
alternatives can be decided by any SWF on two alternatives.

We use our characterization to count the number of SWFs satisfying
IIA, the number of SWFs satisfying IIA and the Weak Pareto Principle,
and the number of SWFs satisfying Arrow's conditions, i.e., IIA and the
Weak Pareto Principle.

Our characterization results in detailed understanding of the SWFs
satisfying Arrow's and Wilson's conditions.  One could argue that all
the functions satisfying these conditions are undesirable voting
systems, as they are very far from any notion of egalitarianism.
However, the richness of the class of such functions as displayed by
their large number may be somewhat surprising.

\subsubsection{Hierarchal Organizations}
A context in which these non-egalitarian functions may be desirable,
or even natural, is that of a hierarchal organization.  Of course
there is no reason to assume organizational decision making should be
transitive (see, e.g., ~\cite{Buchanan54a, Buchanan54b}), let alone
satisfy IIA. However, assuming that {\em individuals} can be rational
then hierarchal organizations {\em can} make rational decisions; for
example, a dictatorship is a simple hierarchal structure that
trivially results in rational organizational decisions, assuming the
dictator is rational. We show that whenever organizational decisions
satisfy Arrow's conditions then the organization is necessarily
hierarchal and characterize the possible structures of
hierarchy. This is not offered as a practical application in the field
of organization design, but rather as a (hopefully) illuminating point
of view on our results.

We claim that given a set of alternatives to be ranked by members of
an organization, an SWF which satisfies Arrow's conditions effectively
determines a ``chain of command'' with respect to each pair of
alternatives $x$ and $y$. Each member of the chain, starting with the
first, can either make a decision (i.e., $x<y$ or $y<x$) or remain
indifferent and pass the decision onto the next in line. A chain can
end with any SWF over two alternatives (e.g., majority vote or a fixed
order), as long as the position of the pair among the rest of the
alternatives has already been fixed. To satisfy Arrow's theorem, these
chains cannot be arbitrary but rather have rich structure. In
particular, there must exist one individual (Arrow's dictator) that is
the first in all chains.

For example, consider a company that must rank alternatives
$\{x,y,z,w\}$. To satisfy Arrow's conditions, it must be that the company
has a CEO that has the prerogative of deciding the ranking. However,
if she is indifferent with respect to $x$ and $y$ then the decision
could be deferred to some subaltern executive, and if he is also
indifferent then $x$ is ranked above $y$. At the same time, it could
also be that when the CEO is indifferent with regards to the ranking
of $z$ and $w$ then a majority vote is taken among the rest of the
employees of the company. This (still not fully specified) SWF meets
all the conditions of Arrow's theorem: IIA, WPP, and transitivity.

The different ``chains of command'' for the different pairs cannot be
constructed arbitrarily. For example, it cannot be the case that the
CEO defers to three distinct individuals the rankings of the three
pairs $(x,y)$, $(y,z)$ and $(x,z)$, since the resulting relation may
not be transitive. The description of the exact conditions of how
these chains may be constructed is the main aim of this article.

\subsubsection{Quantitative Arrovian Theorems}
This complete characterization may also play a role in proving
quantitative versions of Arrow's theorem. Arrow's theorem states that
a non-dictatorial SWF cannot satisfy IIA and transitivity
simultaneously. However, it is natural to ask if such an SWF can have
a transitive outcome for a {\em large proportion} of the voting
profiles. Likewise, a natural question is whether there exists a
non-dictatorial SWF with transitive outcomes, such that the IIA
condition is met for {\em most} voting profiles.

Negative answers to the questions above are given in~\cite{Kalai02,
  Mossel09} where it is shown that if a transitive function satisfies
the IIA property for a large proportion of the profiles then it
necessarily agrees with a dictatorial function on a large proportion
of the profiles. The same holds for an IIA function whose outcome is
transitive for a high proportion of the profiles.

In~\cite{Kalai02} a quantitative Arrow theorem is proved for {\em
  strict preferences}, and only for a restricted sub-class of {\em
  balanced} SWFs (these are functions where for each pair of
alternatives exactly half of the voting profiles result in one
alternative preferred over the other). The proof in~\cite{Mossel09},
which gives a complete quantitative version of Arrow's theorem for
strict preferences, uses a characterization of all the functions
satisfying IIA (for strict preferences). Thus, it proves more generally
that if a function is IIA and close to being transitive then it must
agree on a high proportion of the voting profiles with some functions
that is always IIA and transitive. Hence, it seems that a general
framework which looks at all IIA and transitive functions is crucial
for quantitative proofs of Arrow's theorem.  We expect that the
results proven here will play a role in deriving quantitative versions
of Arrow's theorem, where non-strict preferences are allowed.

\section{Definitions}
Basic notation:
\begin{itemize}
\item $S$ is a finite set of social states (alternatives). In the
  context of this article the interesting regime is $S \geq 3$.
\item $\Pi_S$ is the set of complete transitive binary relations (weak
  preferences) on $S$.

  We depart here from standard notation: given a binary relation $R$,
  we express $(x,y)\in R$ as $x \leq_R y$, rather than the often used
  $xRy$. This gives us the convenience of writing $x <_R y$ and $x =_R
  y$, the meaning of which is clear when $R$ is complete.

  Recall that $R$ is a complete transitive binary relation if it
  satisfies the following conditions:
  \begin{enumerate}
  \item $\forall x,y\in S$, $x\leq_R y$ or $y\leq_R x$ or both. As a
    consequence, $\forall x\in S$ it holds that $x\leq_R x$.
  \item $\forall x,y,z\in S$, if $x\leq_R y$ and $y\leq_R z$ then
    $x\leq_R z$.
  \end{enumerate}

\item $V=\{1,\ldots,N\}$ is a finite set of voters. Although most
  concepts we discuss below (e.g., IIA, Weak Pareto Principle) are
  already non-trivial when $N=1$, our regime of interest is $N \geq 2$.
\item The preference of $i\in V$ is ${\bf R}_i \in \Pi_S$
\item The preferences of all $N$ voters ${\bf R}=({\bf
    R}_{1},\ldots,{\bf R}_{N})$ is called a configuration. Since a
  configuration can also be viewed as a function from the voters to
  preferences, we denote the set of configurations by $\Pi_S^V$.
\item $0$ is the ``null voter,'' who always deems all alternatives
  equal. We introduce the null voter to simplify some of the
  statements to follow.
\item $f:\Pi_S^V\to \Pi_S$ is a social welfare function (SWF).  Note
  that by this definition, all SWFs are transitive, in the sense that
  their image consists of transitive relations.
\end{itemize}

The central property of SWFs that concerns us is Independence of
Irrelevant Alternatives. The following is a standard definition,
introduced by Arrow~\cite{Arrow50}:
\begin{definition}
  An SWF satisfies {\bf Independence of Irrelevant Alternatives (IIA)}
  if, for a given pair $x,y\in S$, it depends only on the voters'
  preferences on that pair. Alternatively, let an SWF $f$ satisfy IIA.
  Then, if configurations ${\bf P} \in \Pi_S^V$ and ${\bf Q} \in
  \Pi_S^V$ agree on $\{x,y\}$, then $f({\bf P})$ agrees with $f({\bf
    Q})$ on $\{x,y\}$.
\end{definition}
Note that any function with only one or two social states trivially
satisfies IIA.

The following definition is also standard.
\begin{definition}
  An SWF satisfies the {\bf Weak Pareto Principle (WPP)} if, whenever all
  the voters agree that $x<y$, then $f$ also agrees on that
  preference: Let $f$ be an SWF satisfying WPP. Then $\forall
  i \in V:\:x<_{R_i}y$ implies $x<_{f({\bf R})}y$.
\end{definition}

Following Wilson, we make the following definition which is a further
weakening of the Weak Pareto Principle:
\begin{definition}
  An SWF satisfies {\bf Citizens' Sovereignty (CS)} if for each pair
  of alternatives $\{x,y\}$, there exists a configuration ${\bf P}$
  such that $x\leq_{f({\bf P})}y$ and there exists a configuration
  ${\bf Q}$ such that $y\leq_{f({\bf Q})}x$.
\end{definition}

Wilson also makes this definition.
\begin{definition}
  An SWF $f$ is {\bf null} if for every configuration ${\bf R}$ and
  pair of alternatives $\{x,y\}$, it holds that $x=_{f({\bf R})}y$.
\end{definition}

Another natural definition is that of a dictatorial function:
\begin{definition}
  An SWF $f$ is {\bf dictatorial} if there exists a dictator in $V$
  whose preferences are a prerequisite for $f$ (perhaps inversely):
  \begin{equation*}
    \begin{array}{c}
      (\exists i\in V)(\forall x,y\in S)(\forall{\bf R}\in \Pi_S^V):\:x\leq _{f({\bf R})}y \to x \leq_{R_i} y\\
      \mbox{or}\\
      (\exists i\in V)(\forall x,y\in S)(\forall{\bf R}\in \Pi_S^V):\:x\leq _{f({\bf R})}y \to y \leq_{R_i} x.
    \end{array}
  \end{equation*}
\end{definition}
Let an SWF $f:\Pi_S^V\to \Pi_S$ be dictatorial, and let $i$ be a
dictator whose preferences are a prerequisite for $f$ (directly,
rather than inversely).  Then, repeating the above, it holds that for
each $x,y$ and ${\bf R}$,
\begin{equation*}
  x \leq _{f({\bf R})} y \to x \leq_{R_i} y.
\end{equation*}
Now, this is equivalent to:
\begin{equation*}
  x <_{R_i} y \to x<_{f({\bf R})}y,
\end{equation*}
so that one could, alternatively, define a dictator as a voter to
whose strict preferences $f$ conforms. The relations among members of
each equivalence class of $R_i$ are then left to be otherwise
determined.

This discussion motivates the following definition:
\begin{definition}
  Let $f:\Pi_S^V\to \Pi_S$ be a dictatorial SWF with dictator $i$ and
  $g:\Pi_S^{V\setminus \{i\}}\to \Pi_S$ to be an SWF, with $i$ removed
  from the set of voters.

  Given a configuration ${\bf R}\in \Pi_S^V$, we will write ${\bf
    R}_{-i}$ for the element of $\Pi_S^{V\setminus \{i\}}$ which
  agrees with ${\bf R}_{-i}$ for all voters in $V\setminus \{i\}$.

  We say that $f$ {\bf defers to $g$ over $i$} if, whenever $x=_{R_i}
  y$, $f({\bf R})$ agrees with $g({\bf R}_{-i})$ over $\{x,y\}$:
\begin{equation*}
  (\forall x,y \in S)(\forall {\bf R}\in \Pi_S^V):\: x=_{R_i}y \to
  \left( x\leq_{f({\bf R})}y \:\leftrightarrow\: x\leq_{g({\bf R}_{-i})}y \right)
\end{equation*}
\end{definition}
Intuitively, the phrase ``$f$ defers to $g$ over $i$'' means that the
function $g$ is used to rank any pair for which the dictator $i$ is
undecided. It is easy to see that for any dictatorial $f$ with
dictator $i$ there may only exist a single $g$ such that ``$f$ defers
to $g$ over $i$''.

Let a ``cleric'' be characterized by some {\em fixed} opinion $C\in
\Pi_S$, which is independent of all the voters' preferences.
\begin{definition}
  An SWF $f$ is {\bf clerical with cleric $C$} if $C$'s (fixed)
  preferences are a prerequisite for $f$:
\begin{equation*}
  (\forall x,y\in S)(\forall {\bf R}\in \Pi_S^V):\:x\leq _{f({\bf R})}y\to x\leq_Cy.
\end{equation*}
\end{definition}
Alternatively, if $C$ is a cleric for $f$, then $C$'s strict
preferences are followed by $f$, for all configurations. Note that for
$f$ clerical with cleric $C$, if $\{x,y\}$ are such that
$x\leq_{f({\bf R})}y$ for some configurations, and $y\leq_{f({\bf
    R})}x$ for some configurations (so that there is no strict fixed
preference between them), then $x=_Cy$, and $x$ and $y$ belong to the
same equivalence class of $C$. Note also that every SWF is clerical
with the trivial cleric who deems all alternatives equal.

Let $A \subset S$ be a set of alternatives and $f$ an SWF. The
function $f$ restricted to $A$, denoted $f|_A$, is the function $f|_A
: \Pi_A^V \to \Pi_A$ defined as follows. Given ${\bf R} = ({\bf
  R}_1,\ldots,{\bf R}_n) \in \Pi_A^V$, let ${\bf R'} = ({\bf
  R}'_1,\ldots,{\bf R}'_n) \in \Pi_S^V$ agree with them on all the
pairwise preferences in $A$.  Then for $x,y \in A$ we have $x \leq
_{f|_A({\bf R})} y$ iff $x \leq _{f({\bf R'})} y$. Note that if (but
not only if) IIA holds then $f|_A$ is well defined, i.e. it does not
depend on the choice of ${\bf R'}$. This point has been hitherto noted
by Blau~\cite{Blau71}.

\begin{definition}
  An SWF $f:\Pi_S^V\to \Pi_S$ is {\bf clerical-dictatorial} if it is
  null, if $S$ is of size at most two, or if there exists a cleric $C$
  such that all of the following hold:
  \begin{enumerate}
  \item $f$ is a clerical function with cleric $C$.
  \item Let $A$ be any of $C$'s equivalence classes.
    \begin{enumerate}
    \item $f|_A$ is well defined, and is either null, dictatorial, or
      of size $2$.
    \item If $f|_A$ is dictatorial then it defers to a
      {\bf clerical-dictatorial} function.
    \end{enumerate}
  \end{enumerate}
\end{definition}

\section{Our Theorem}
\label{section:theorem}
In terms of the definitions above, Arrow's theorem~\cite{Arrow50} is
the following:
\begin{theorem}[Arrow]
  Any SWF which satisfies IIA and WPP is dictatorial (with the
  dictator's decrees followed verbatim rather than in reverse).
\end{theorem}

Wilson's first theorem (theorem 3 in~\cite{Wilson72}) is the
following:
\begin{theorem}[Wilson]
  Any SWF which satisfies IIA and CS is either dictatorial or null.
\end{theorem}

Wilson's second theorem (theorem 5 in~\cite{Wilson72}) is the following:
\begin{theorem}[Wilson]
  \label{thm:wilson}
  Any SWF which satisfies IIA is clerical with some cleric $C$, so
  that restricted to $C$'s equivalence classes of size $>2$, it is
  either dictatorial or null.
\end{theorem}

Note that none of the above is a complete characterization of the
functions satisfying the respective conditions.
\begin{example}
  Consider the following function $f$ on three social alternatives,
  $x$, $y$ and $z$, with two voters. Let voter 1 be a dictator, so
  that his strict preferences are followed verbatim. When voter 1
  deems all three alternatives equal, then
\begin{itemize}
\item If voter 2 also votes $x=y=z$ then $f$ outputs $x=y=z$.
\item Otherwise (for example, when voter 2 votes $x<y=z$), $f$ outputs
  $x=y<z$.
\end{itemize}
This function, although it is dictatorial (and even satisfies WPP),
clearly does not satisfy IIA.
\end{example}

By the definitions above, any SWF is transitive, in the sense that the
relations in its image are transitive. If one allows for
non-transitive preferences in the image, then it possible to construct
a dictatorial function that is not transitive.
\begin{example}
  Let $f$ be an SWF on three social alternative and two voters, with
  the first voter taking the r\^ole of dictator. Then $f$ follows
  voter 1's strict preferences. When 1 deems two alternatives equal,
  then $f$ does too, if voter 2 agrees. When voter 2 disagrees, then,
  depending on the two alternatives in questions, $f$ outputs $x<y$ if
  the two alternatives are $x,y$, $y<z$ if the two alternatives are
  $y,z$ and $z<x$ if the two alternatives are $x,z$.  Thus, when voter
  1 votes $x=y=z$ and voter 2 votes $x<y<z$, then $f$ outputs $x<y<z$,
  but also $z<x$, and hence is not transitive.  Note that this
  function satisfies IIA and WPP.
\end{example}

We provide a complete characterization of SWFs satisfying the
conditions in theorems of Wilson and of Arrow.

Our ``if and only if'' version of Arrow's theorem is:
\begin{theorem} \label{thm:main_arrow} An SWF satisfies IIA and WPP if
  and only if is dictatorial (with the dictator's decrees followed
  verbatim rather than in reverse), and defers to a
  clerical-dictatorial function.
\end{theorem}

Our ``if and only if'' version of Wilson's first theorem is:
\begin{theorem} \label{thm:main_wilson1} An SWF satisfies IIA and CS
  if and only if it is null or dictatorial, and in the latter case
  defers to a clerical-dictatorial function.
\end{theorem}

Our ``if and only if'' version of Wilson's second theorem is:
\begin{theorem} \label{thm:main} An SWF satisfies IIA if and only if
  it is clerical-dictatorial.
\end{theorem}

We prove~\eqref{thm:main} below, and show that theorems~\eqref{thm:main_arrow}
and~\eqref{thm:main_wilson1} follow as corollaries.

\subsection{IIA $\implies$ Clerical-Dictatorial}

In this section we prove the following theorem:

\begin{theorem} \label{thm:IIA2CD}
Any SWF which satisfies IIA is clerical-dictatorial.
\end{theorem}

A key concept in the proof will be the relation $C(f)$, defined as follows:
$C(f)$  is the relation signifying that,
between social states $x$ and $y$, for some configuration,
$x\leq_fy$. Formally,
\begin{equation*}
  x\leq_{C(f)}y\quad\leftrightarrow\quad\exists {\bf R}\in \Pi_S^V:\:
x\leq_{f({\bf R})}y
\end{equation*}

Note that $x<_{C(f)}y$ iff for all configurations, $x<_{f({\bf R})}y$.

\begin{figure}[h]
  \label{fig:clerical}
  \centering
  \includegraphics{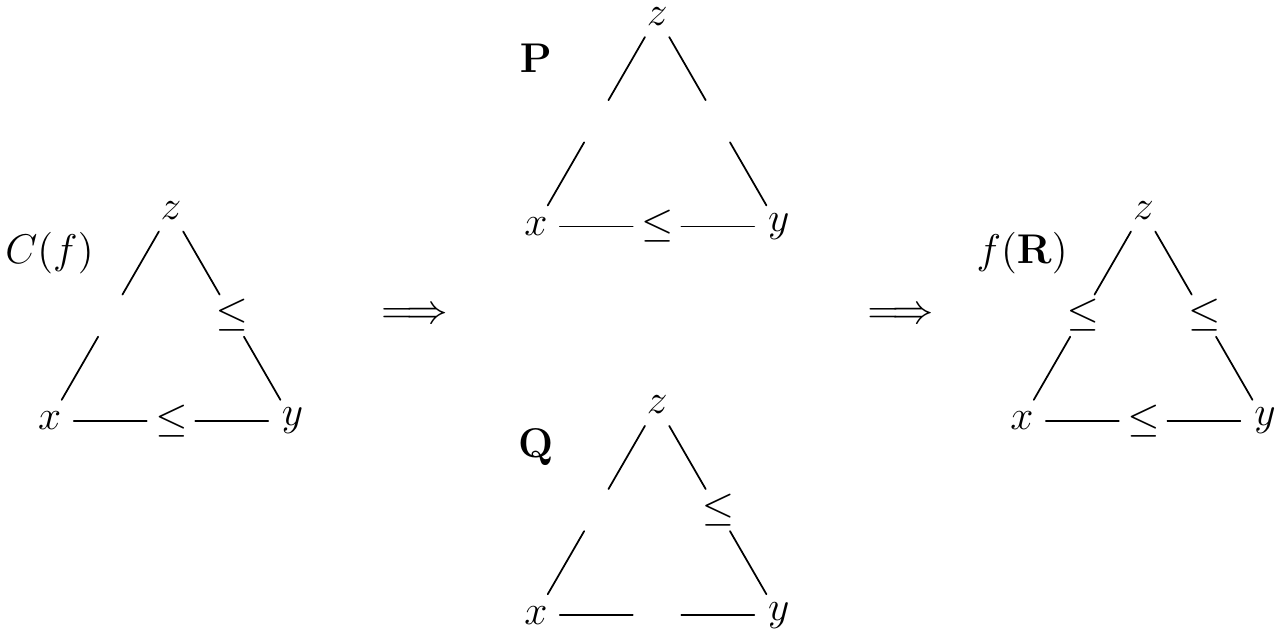}
  \caption{Transitivity of $C(f)$.}
\end{figure}

We will use the following lemma which is essentially due to Wilson.
\begin{lemma}
  \label{lemma:clerical}
  $C(f)$ is complete and transitive.
\end{lemma}
\begin{proof}
  Completeness: let ${\bf R}$ be some configuration. Then either
$x\leq_{f({\bf R})}y$ or $y\leq_{f({\bf R})}x$. Hence $x\leq_{C(f)}y$ or
$y\leq_{C(f)}x$, and $C(f)$ is complete.

Transitivity (see Fig.~\ref{fig:clerical}): let $x\leq_{C(f)}y$ and
$y\leq_{C(f)}z$.  Then there exist configurations ${\bf P}$ and ${\bf
  Q}$ s.t. $x\leq_{f({\bf P})}y$ and $y\leq_{f({\bf Q})}z$. Let ${\bf
  R}$ be a configuration that agrees with ${\bf P}$ on $\{x,y\}$ and
with ${\bf Q}$ on $\{y,z\}$. Such a configuration always exists, since
preferences on $\{x,y\}$ and $\{y,z\}$ can always be completed to a
transitive relation by some ordering of $\{x,z\}$.

Because of IIA, $x\leq_{f({\bf R})}y$ and $y\leq_{f({\bf R})}z$.
Finally, because $f({\bf R})$ is transitive, it follows that
$x\leq_{f({\bf R})}z$, and hence, by the definition of $C(f)$,
$x\leq_{C(f)}z$.
\end{proof}

Note that since $C(f)$ is a weak preference,
lemma~\ref{lemma:clerical} implies that any SWF that satisfies IIA is
clerical, with cleric $C(f)$.  We now define $I_f(x)$ to be the
equivalence class of $x$ under $C(f)$. By this definition, if $y$ is
in $I_f(x)$, then for some configurations $x\leq_{f}y$, and for some
$y\leq_{f}x$.  Alternatively, there is no fixed strict preference
between $x$ and $y$.

Wilson's result is that $f$, restricted to $I_f(x)$ of size $>2$ is
either null or dictatorial. We now assign each social state $x$ a
dictator $D(f,x)$ as follows.  If $f$ restricted to $I_f(x)$ is
dictatorial, then we let $D(f,x)$ be that dictator.  Otherwise we let
$D(f,x) = 0$ be the null dictator.

The following definition is the key to understanding the structure of
SWFs that satisfy IIA:
\begin{definition}
  Given an SWF $f$ and social state $x$, we define a series of $N+1$
  SWFs, $f_0^x$ to $f_N^x$, recursively:
  \begin{itemize}
  \item The basis of the recursion is trivial:
    \[
    f_0^x:=f, \quad V_0 := V, \quad S_0 := S.
    \]
  \item For $i\geq 0$:

    Denote $w = D(f_i,x)$ and $A=I_{f_i^x}(x)$.
    \begin{itemize}
    \item $S_{i+1}:=A$.
    \item If $w=0$ then $V_{i+1}:=V_i$ and $f_{i+1}^x:=f_i^x|_A$.
    \item If $w\neq 0$ then $V_{i+1}:=V_i\setminus\{w\}$ and
      $f_{i+1}^x$ is the function that $f_i^x|_A$ defers to over $w$.
    \end{itemize}
  \end{itemize}

\end{definition}
\begin{lemma}
  \label{lemma:IIA_decomposition}
  When $f$ satisfies IIA then $f_i^x$ is well defined and satisfies
  IIA for $0 \leq i \leq N$.
\end{lemma}
\begin{proof}
  Clearly $f_0^x$ is well defined and satisfies IIA when $f$
  does. Assume that $f_i^x$ is well defined and satisfies IIA. We will
  show that (1) $f_{i+1}^x$ is well defined and (2) satisfies IIA.

  \begin{enumerate}
  \item Since $f_i^x$ satisfies IIA, then $f_i^x$'s preferences on the
    part of its domain that it {\em doesn't} share with $f_{i+1}^x$,
    $S_i\setminus S_{i+1}$, do not influence $f_i^x({\bf R}_i)$ for
    states in $S_{i+1}$, and $f_{i+1}^x$ is well defined.
  \item Since $f_{i+1}^x$ equals to $f_i^x$ on a restricted domain, it
    too satisfies IIA.
  \end{enumerate}
\end{proof}

To better understand the functions $f_i^x$, note that
\begin{enumerate}
\item $S_{i+1}\subseteq S_i$.
\item $V_{i+1}\subseteq V_i$ and $|V_i|-1\leq |V_{i+1}| \leq|V_i|$.
\item If $y$ is an element of $I_{f_i^x}$ then
  $f_{i+1}^y=f_{i+1}^x$. Alternatively, if $x$ and $y$ are in $S_i$,
  then $f_i^x=f_i^y$.
\end{enumerate}

Noting that if $y\in S_i$ then $f_i^x=f_i^y$, we may state the
following immediate lemma:
\begin{lemma}
  \label{lemma:defers}
  If $f_i^x$, restricted to some equivalence class $I_{f_i^x}(y)$, is
  dictatorial, then it defers to $f_{i+1}^y$.
\end{lemma}

Some insight may be gained by an ``algorithmic'' view of how $f$
``decides'' on a pair $\{x,y\}$.  Starting with $i=0$, this is
performed as outlined in Algorithm 1.
\begin{algorithm}
\label{algorithm}
\caption{$f_i^x(x,y)$}
\begin{algorithmic}[htb]
\IF {the cleric $C(f_i^x)$ has a strict preference for $\{x,y\}$}
  \STATE  $f$ conforms to it
\ELSE
  \STATE  {It must be the case that $y\in I_{f_i^x}$. Let $w = D(f_i,x)$.}
  \IF{$w=0$}
    \STATE{$f_i^x(x,y)$ is null, or else $S_i=\{x,y\}$ and $f_i^x(x,y)$
         is a general function of the preferences on $x,y$
         of the voters in $V_i$.}
  \ELSIF{$w \neq 0$ and $w$ has strict preference on $x,y$}
    \STATE  $f_i^x$ follows $w$
  \ELSE
     \STATE {$x=_{R_w}y$ and the decision is deferred to $f_{i+1}^x$}
  \ENDIF
\ENDIF
\end{algorithmic}
\end{algorithm}

\begin{lemma}
  \label{lemma:almost_constant}
  $f_N^x$ is null or has at most two social states.
\end{lemma}
\begin{proof}
  Note that if for some $i$ we have $V_{i+1}=V_i$ then $f_N^x=f_i$ has
  at most two social states or is null.  Since $V_{i+1}$ can be
  smaller than $V_i$ for no more than $N$ values of $i$, the proof
  follows.
\end{proof}
Therefore trivially $f_N^x$ is by definition {\bf
  clerical-dictatorial}.  We may now prove theorem~\ref{thm:IIA2CD}:
any SWF which satisfies IIA is clerical-dictatorial.

\begin{proof}
  Let $f$ satisfy IIA. Then the function $f_i^x$ is well defined for
  all $x\in S$ and $0\leq i \leq N$ by
  lemma~\ref{lemma:IIA_decomposition}. Furthermore, it satisfies IIA.

  Now, for all $x$, $f_N^x$ is clerical-dictatorial by
  lemma~\ref{lemma:almost_constant}. Assume that $f_{i+1}^x$ is
  clerical dictatorial for all $x$. We will show that $f_i^x$ is
  clerical-dictatorial, and thus prove by induction that $f=f_0^x$ is
  too.

  If $f_i^x$ is null or has only one or two social states then it is
  clerical dictatorial. Otherwise, it is clerical with cleric
  $C(f_i^x)$. Let $A=I_{f_i^x}(y)$ be some equivalence class of
  $C(f_i^x)$ of size $>2$. Then
  \begin{enumerate}
  \item Since $f_i^x$ satisfies IIA then $f_i^x|_A$ is well
    defined. As a consequence of the same fact, by
    theorem~\ref{thm:wilson} (Wilson), $f_i^x|_A$ is either
    dictatorial, null, or of size $2$.
  \item If $f_i^x|_A$ is dictatorial then it defers to $f_{i+1}^y$ by
    lemma~\ref{lemma:defers}. By the inductive hypothesis $f_{i+1}^y$
    is clerical-dictatorial.
  \end{enumerate}
  Thus $f_i^x$ satisfies the conditions for being clerical-dictatorial.

\end{proof}

\subsection{Clerical Dictatorial $\implies$ IIA}

\begin{theorem}
A clerical-dictatorial SWF satisfies IIA.
\end{theorem}
\begin{proof}
  The proof is by induction on the number of voters.  The base case of
  no voters is trivial.  Now assume that there are $n \geq 1$ voters
  and let $f:\Pi_S^V\to \Pi_S$ be a clerical-dictatorial SWF, with
  cleric $C$. Let ${\bf P}$ and ${\bf Q}$ be configurations that agree
  on $\{x,y\}$.  We would like to show that $f({\bf P})$ agrees with
  $f({\bf Q})$ on $\{x,y\}$.

  If $C(f)$ is strict about $\{x,y\}$ then $f$ is constant with regard
  to them, so in particular $f({\bf P})$ agrees with $f({\bf Q})$ on
  $\{x,y\}$.  Otherwise, $y\in I_f(x)$, and we denote $A=I_f(x)$. By
  the definition of clerical-dictatorial functions, $f|_A$ is well
  defined, so that the ranking of $x$ and $y$ depends only on the
  voters' opinions with regards to states in $A$.

  The case where $f|_A$ is null is trivial.  Likewise, if $I_f(x)$ is
  of size two then we're done as the ranking between $x$ and $y$
  depends only on the individual preferences between them.

  In the remaining case $f|_A$ is dictatorial. Let $i$ be the dictator
  and $g$ the clerical-dictatorial function that it defers to.  Now if
  the dictator $i$ has strict ranking between $x$ and $y$ in $P_i$ and
  $Q_i$ then $f({\bf P})$ and $f({\bf Q})$ follow the dictator and
  agree on $\{x,y\}$.  Otherwise, the relative ranking of $x$ and $y$
  is determined by $g$, which is IIA by the inductive hypothesis, and
  so $f({\bf P})$ and $f({\bf Q})$ agree on $\{x,y\}$ as needed.

\end{proof}

\subsection{Arrow's Model and Wilson's Model}
The following two theorems, also mentioned above, are simple applications of our
main result to the domains of Arrow's and Wilson's models, i.e. the set
of SWFs satisfying, besides IIA, WPP and CS, respectively.

Our ``if and only if'' version of Arrow's theorem is:
\begin{theorem*}[\ref{thm:main_arrow}]
  An SWF satisfies IIA and WPP if and only if is dictatorial (with the
  dictator's decrees followed verbatim rather than in reverse), and
  defers to a clerical-dictatorial function.
\end{theorem*}
\begin{proof}
First direction: an SWF which satisfies IIA and WPP must be dictatorial (with
the dictator's decrees followed verbatim rather than in reverse), by Arrow's
theorem. By our theorem, it must be clerical-dictatorial and hence
defer to a clerical dictatorial function.

Second direction: an SWF which is dictatorial (again, verbatim) certainly
satisfies WPP. If it is also clerical-dictatorial, then it also satisfies IIA.
\end{proof}

Our ``if and only if'' version of Wilson's first theorem is:
\begin{theorem*}[\ref{thm:main_wilson1}]
  An SWF satisfies IIA and CS if and only if is null or dictatorial,
  and in the latter case defers to a clerical-dictatorial function.
\end{theorem*}
The proof of this theorem is essentially identical to the previous.
\begin{proof}
First direction: an SWF which satisfies IIA and CS must be null or dictatorial,
by Wilson's theorem. By our theorem, it must be clerical-dictatorial and hence
defer to a clerical dictatorial function.

Second direction: an SWF which is dictatorial or null certainly
satisfies CS. If it is also clerical-dictatorial, then it also satisfies IIA.
\end{proof}

\section{The Number of Social Welfare Functions Satisfying IIA}
The characterization of (transitive) social welfare functions satisfying
IIA as clerical-dictatorial facilitates counting the number of such
functions, for low numbers of social states.

The sequence starting with $1,1,13,75,541,\ldots$, which counts the number of
clerics (weak preferences) over $S$ social states, is known as the
{\em Ordered Bell Numbers} or {\em Fubini Numbers} \cite{wiki:StrictWeakOrder}.
There is no known simple expression for this sequence, recursive or otherwise.

Denote by $q_s(v)$ the number of SWFs satisfying IIA, over $s$
social states and $v$ voters. Then $q_s(0)$ are the Ordered Bell Numbers.
For small $s$, it is possible to form a recursive expression for $q_s(v)$. We do
that by first choosing a cleric $C$, and then counting the number
of functions for each equivalence class.

For $s=1$, only one function is possible. For $s=2$, any function is
transitive and satisfies IIA, and therefore
\begin{equation*}
  q_2(v)=3^{3^v}.
\end{equation*}

For $s=3$
\begin{equation*}
  q_3(v) = 6 + 6 \cdot (q_2(v) - 2) + (1 + 2vq_3(v-1)),
\end{equation*}
with $q_3(0)=13$.
The first addend counts the number of strict clerics, that
is, the number of fixed strict functions. There are six such clerics,
and each equivalence class is of size one and therefore has only one
possible function.

The second addend counts the number
of functions for which $C$ has an equivalence class of size one and an
equivalence class of size two. There are six such clerics, and for each one
there are $q_2(v)$ possible functions for the size two equivalence class. Of
these, we have to subtract the two constant strict ones which we counted in
the third addend.

The third addend counts the number of functions possible when $C$ is
the unique cleric that deems all three social states equal. Then,
either the function is constant and equal to $C$, or else it is
dictatorial. There are $v$ possible dictators, and each one can be
followed directly or inversely. Once that is chosen, there are
$q_3(v-1)$ options for the function that that dictator defers to.

Note that $q_3(v)$ is dominated by the second addend, or by the number of
functions over two social states.

A similar analysis, although increasingly complicated, can be done
for larger $s$. For example, for $s=4$,
\begin{equation*}
  q_4(v)=24+36\cdot(q_2(v)-2)+6\cdot(q_2(v)-2)^2+8\cdot(1+2vq_3(v-1))+
  (2vq_4(v-1)+1).
\end{equation*}
Here, $q_4(v)$ is dominated by the third addend, or the number of functions
for which $C$ has two equivalence classes of size two.

\subsection{Arrow's Model}
The calculation of the number of SWFs satisfying Arrow's conditions is simple,
when expressed in terms of $q_s(v)$.
Let $r_s(v)$ be the number of SWFs satisfying IIA and the
Weak Pareto Principle. Then $r_s(0)$ is equal to $q_s(0)$, or the Ordered Bell
Numbers.

For $v>1$ and $s > 1$:
\begin{equation*}
  r_s(v)=vq_s(v-1).
\end{equation*}
Since we must pick a dictator and the clerical-dictatorial function that is
deferred to.

\subsection{Wilson's Model}
The calculation in this case is also simple, when expressed in terms
of the $q_s(v)$.  Let $p_s(v)$ be the number of SWFs satisfying IIA
and Citizens' Sovereignty. Then $p_s(0)$ is one, since the only
possible function is the null function.

For $s=2$, any function is admissible, except the two that have constant
strict rankings:
\begin{equation*}
  p_2(v)=3^{3^v}-2.
\end{equation*}

For $s=3$
\begin{equation*}
  p_3(v) = 1 + 2vq_3(v-1)
\end{equation*}
which is the last term in the expression for $q_3(v)$, the term which
counts the number of functions in the case where $C(f)$ deems all
alternatives equal.  As mentioned above, this conditions is equivalent
to Citizens' Sovereignty. Likewise, for any $s>2$:
\begin{equation*}
  p_s(v) = 1 + 2vq_s(v-1).
\end{equation*}

\bibliographystyle{abbrv} \bibliography{scf}

\begin{thebibliography}{1}

\bibitem{Arrow50}
K.~Arrow.
\newblock A difficulty in the concept of social welfare.
\newblock {\em J. of Political Economy}, 58:328--346, 1950.

\bibitem{Arrow63}
K.~Arrow.
\newblock {\em Social choice and individual values}.
\newblock John Wiley and Sons, 1963.

\bibitem{Blau71}
J.~Blau.
\newblock Arrow's theorem with weak independence.
\newblock {\em Econometrica}, 538(152):413--420, 1971.

\bibitem{Buchanan54a}
J.~Buchanan.
\newblock {Individual choice in voting and the market}.
\newblock {\em The Journal of Political Economy}, 62(4):334--343, 1954.

\bibitem{Buchanan54b}
J.~Buchanan.
\newblock {Social choice, democracy, and free markets}.
\newblock {\em The Journal of Political Economy}, 62(2):114--123, 1954.

\bibitem{Kalai02}
G.~Kalai.
\newblock {A Fourier-theoretic perspective on the Concordet paradox and Arrow's
  theorem}.
\newblock {\em {Adv.\ in Appl.\ Math.}}, 29(3):412--426, 2002.

\bibitem{Mossel09}
E.~Mossel.
\newblock A quantitative arrow theorem.
\newblock Available at the Arxiv 0903.2574, 2009.

\bibitem{wiki:StrictWeakOrder}
Wikipedia.
\newblock Strict weak ordering --- {W}ikipedia{,} the free encyclopedia, 2009.
\newblock [Online; accessed 15-July-2009].

\bibitem{Wilson72}
R.~Wilson.
\newblock Social choice theory without the pareto principle.
\newblock {\em Journal of Economic Theory}, 5(3):478--486, December 1972.

\end{thebibliography}

\end{document}